\documentclass{article}

\usepackage{amsmath, amscd, amsfonts, amssymb, amsthm}

\newtheorem{thm}{Theorem}[section]
\newtheorem{cor}[thm]{Corollary}
\newtheorem{lem}[thm]{Lemma}
\newtheorem{prop}[thm]{Proposition}
\begin{document}
\title{The Noisy Voter-Exclusion Process\footnote{Research supported in part by NSF grant DMS-00-70465}}
\author{Paul Jung}
\maketitle
\begin{abstract}
The symmetric exclusion process and the voter model are two
interacting particle systems for which a dual finite particle
system allows one to characterize its invariant measures. Adding
spontaneous births and deaths to the two processes still allows
one to use the dual process to obtain information concerning the
original process. This paper introduces the noisy voter-exclusion
process which generalizes these processes by allowing for all of
these interactions to take place. The dual process is used to
characterize its invariant measures under various circumstances.
Finally, an ergodic theorem for a related process is proved using
the coupling method.
\end{abstract}
\section{Introduction}
The voter model is an interacting particle system introduced
independently by Clifford and Sudbury(1973) and Holley and
Liggett(1975).  In particular it is a spin system (see
Liggett(1985)) with rates given by
\begin{equation*}
c(x,\eta)=\left\{
\begin{array}{ll}
\sum_y q_v(x,y)\eta(y)&\text{if }\eta(x)=0,\\
\sum_y q_v(x,y)[1-\eta(y)]&\text{if }\eta(x)=1,
\end{array}
\right.
\end{equation*}
where $q_v(x,y)\ge 0$  and $\sup_x \sum_y q_v(x,y)<\infty$ for
$x,y\in \mathcal{S}$.

To describe the voter model in a more intuitive manner let
$\mathcal{S}$ be a countable set for which a voter resides at each
site in the set.  The voter at site $x$ waits an exponential time
with mean $[\sum_y q_v(x,y)]^{-1}$ at which point it chooses one
of its neighbors with probability $q_v(x,y)/\sum_z q_v(x,z)$ and
subsequently takes the opinion (either $1$ or $0$) of $y$.

Schwartz(1976) introduced the $\beta$-$\delta$ process, a particle
system which modifies the well known symmetric exclusion process
with transition rates $q_e(x,y)$, by allowing a birth with
exponential rate $\beta(x)$ when there is a $0$ at site $x$, and a
death with rate $\delta(x)$ when there is a $1$ at site $x$.

Define the transition rates $q(x,y)=q_e(x,y)+q_v(x,y)$, and let
$q_x=\sum_y q(x,y)$. We combine the voter model and the
$\beta$-$\delta$ process to obtain a new process which much
satisfy the following: (a) $\mathcal{S}$ is irreducible with
respect to $q(x,y)$, (b) $q_e(x,y)=q_e(y,x)$, (c) $\sup_x
q_x<\infty$, and (d) $\inf_x q_x>0$. Add to this the transition
rates $\beta(x)$ and $\delta(x)$ where $\sup_x (\beta(x)
+\delta(x))<\infty$. Condition (d) is not necessary, but it is
convenient for the purposes of our discussion. We will call such a
process a noisy voter-exclusion process (NVE process). The NVE
process is a particular example of a spin system with stirring,
also known as a reaction-diffusion process.  In the physics
literature these processes are known as having Glauber-Kawasaki
dynamics.

In the setting of the NVE process, the voter at $x$ waits an
exponential time with mean $q_x$ at which point it again chooses a
neighbor with probability $q(x,y)/q_x$, but now the voter decides
to either switch places with $y$ with probability
$q_e(x,y)/[q_e(x,y)+q_v(x,y)]$ or, as before, take the opinion of
$y$ with probability $q_v(x,y)/[q_e(x,y)+q_v(x,y)]$.  In addition
to this, a voter at $x$ with opinion $0$ decides to spontaneously
change its opinion to $1$ with exponential rate $\beta(x)$, and a
voter at $x$ with opinion $1$ spontaneously changes its opinion to
$0$ with rate $\delta(x)$.

Let
\begin{equation*}
\eta_x(u)=\left\{
\begin{array}{ll}
\eta(u)&\text{if }u\neq x\\
1-\eta(u)&\text{if }u=x
\end{array}\right.
\end{equation*}
and
\begin{equation*}
\eta_{xy}(u)=\left\{
\begin{array}{ll}
\eta(y)&\text{if }u=x\\
\eta(x)&\text{if }u=y\\
\eta(u)&\text{if }u\neq x,y.
\end{array}\right.
\end{equation*}
Using the results of Chapter I in IPS (Liggett(1985)), the
generator for an NVE process is given by the closure of the
following operator on $\mathcal{D}$, the set of all functions on
$\{0,1\}^\mathcal{S}$ that depend on finitely many coordinates:
\begin{equation*}
{\Omega}f(\eta)=\sum_{\eta(x)=1,\eta(y)=0}{q_e(x,y)[f(\eta_{xy})-f(\eta)]}
+\sum_{x}{c(x,\eta)[f(\eta_{x})-f(\eta)]}
\end{equation*}
where
\begin{equation*}
c(x,\eta)=\left\{
\begin{array}{ll}
\beta(x)+\sum_y q_v(x,y)\eta(y)&\text{if }\eta(x)=0,\\
\delta(x)+\sum_y q_v(x,y)[1-\eta(y)]&\text{if }\eta(x)=1.
\end{array}
\right.
\end{equation*}
We will call the corresponding semigroup $S(t)$.

If $\beta(x)=\delta(x)\equiv 0$ then we will say that we have a
voter-exclusion process. One may also refer to this as the voter
model with stirring.   A previous study (Belitsky \textit{et
al.}(2001)) has been done concerning the ergodic theory of the
voter-exclusion process in the case where $\mathcal{S}=\mathbb{Z}$
and $q_e(x,y)$ is not necessarily symmetric, but there is no
overlap with the results of this paper.

If $q_e(x,y)\equiv 0$ then we just get the noisy voter model.
Granovsky and Madras(1995) study some important equilibrium
functionals and critical values of the noisy voter model, but only
for the case where $\beta$ and $\delta$ are constant. We, on the
other hand, will study the invariant measures of the NVE process
where $\beta(x)$ and $\delta(x)$ are in general not constant.

In Chapters V and VIII of IPS, one can find a complete
characterization of the extremal invariant measures and their
domains of attraction for the voter model (Holley and
Liggett(1975)) and the symmetric exclusion process (Liggett(1973)
and Spitzer(1974)) respectively. Schwartz(1976) does the same for
the $\beta$-$\delta$ process. These results are all based upon the
existence of a certain dual finite particle process
(Spitzer(1970)) and a certain monotonicity concerning this dual
process.  In particular, $S(t)\hat{\nu}_\alpha(A)$ defined below
is nondecreasing in $t$ for the voter model and nonincreasing in
$t$ for the symmetric exclusion process. For the NVE process, a
dual still exists, however, there is no monotonicity concerning
the dual so we will have to use other techniques in order to
classify the invariant measures under various conditions. Assume
throughout that $q_v(x,y)>0$ for some $x,y\in\mathcal{S}$ since
all other cases have been studied by Schwartz(1976).

We start with some definitions. Let $\mathcal{P}$ denote the set
of probability measures on $X=\{0,1\}^\mathcal{S}$. The set
$\mathcal{I}$ will denote the invariant measures for a given NVE
process, and $\mathcal{I}_e$ will be its extreme points.

If we denote the set of nonnegative harmonic functions bounded by
$1$ on $\mathcal{S}$ as
\begin{equation*}
\mathcal{H}=\left\{\alpha:\mathcal{S}\rightarrow [0,1] \text{ such
that } \sum_y q(x,y)\alpha(y)=q_x\alpha(x) \text{ for all
}x\right\},
\end{equation*}
then we can define $\nu_\alpha$ to be the product measure on $X$
with marginals $\nu_\alpha\{\eta:\eta(x)=1\}=\alpha(x)$. Let
$\mu_\alpha=\lim_{t\rightarrow\infty}\nu_\alpha S(t)$. Theorem
\ref{thm3.4} below will show that these limits exist.

Let $\mathcal{S}_n=\mathcal{S}^n\backslash\{\vec{x}:x_i=x_j\text{
for some }i<j \}$. If $E_t=(x_t,y_t)\in\mathcal{S}_2$ is the
finite, two particle exclusion process with transition rates
$q(x,y)$ then define the functions $q_v$ and $q_e$ on
$\mathcal{S}_2$ by $q_v(E_t)=q_v(x_t,y_t)+q_v(y_t,x_t)$ and
$q_e(E_t)=q_e(x_t,y_t)+q_e(y_t,x_t)=2q_e(x_t,y_t)$.

Suppose $X(t)$ and $Y(t)$ are independent continuous time Markov
chains on $\mathcal{S}$ with transition rates $q(x,y)$ and denote
$p_t(x,y)=P^x(X(t)=y)$. Let $\Lambda=\{\omega| \int_0^\infty
\beta(X(t))+\delta(X(t))\,dt<\infty\}$. For
$\alpha\in\mathcal{H}$, $\alpha(X(t))$ is a bounded martingale so
$\lim_{t\rightarrow\infty}\alpha(X(t))$ exists with probability
one. We can define an equivalence relation $R$ on $\mathcal{H}$ by
\begin{equation*}
\alpha_1R\alpha_2\text{ if
}\lim_{t\rightarrow\infty}[\alpha_1(X(t))-\alpha_2(X(t))]=0 \text{
almost surely on }\Lambda.
\end{equation*}
$\mathcal{H}_R$ is any set of representatives of the equivalence
classes determined by $R$.

Let $\mathcal{E}$ be the following event:
\begin{equation*}
\{\text{there exists }t_n\rightarrow\infty\text{ such that
}X(t_n)=Y(t_n)\}.
\end{equation*}
Then we will say that $\mathcal{H}^*$ is the set of all
$\alpha\in\mathcal{H}$ such that
\begin{equation*}
P^{\{x,y\}}(\lim_{t\rightarrow\infty} \alpha(X(t))=0\text{ or
}1\text{ on }\mathcal{E})=1\text{ for all }x,y\in\mathcal{S},
\end{equation*}
and $\mathcal{H}^*_R$ is again the set of equivalence classes on
$\mathcal{H}^*$.

Define the following function on $\mathcal{S}^2$,
\begin{equation*}
g(x,y)=P^{(x,y)}[X(t)=Y(t)\text{ for some }t>0].
\end{equation*}
Note that if $g(x,y)=1$ for some $(x,y)\in\mathcal{S}_2$ then by
irreducibility $g(x,y)\equiv 1$ (For more detail concerning this
see Lemma VIII.1.18 in IPS).

We are now in a position to state the theorems:
\begin{thm}\label{thm3.0}
An NVE process is ergodic if and only if
\begin{equation}\label{eq0}
P^x[\int_0^\infty \beta(X(t))+\delta(X(t))\, dt=\infty]=1 \text{
for all }x\in\mathcal{S}.
\end{equation}
\end{thm}

\begin{thm}\label{thm3.1}
Suppose $\mu\in\mathcal{P}$ and $\delta_0, \delta_1$ are the point
masses on all $0$'s and all $1$'s. Assume that (\ref{eq0}) does
not hold and that
\begin{equation}\label{schwartzcond}
P^E[\int_0^\infty q_v(E_t)\, dt=\infty]=1 \text{ for all
}E\in\mathcal{S}_2.
\end{equation}
Then

(a) $\lim_{t\rightarrow\infty} \delta_0S(t)=\mu^0$ and
$\lim_{t\rightarrow\infty} \delta_1S(t)=\mu^1$ exist,

(b) $\mathcal{I}_e=\{\mu^0,\mu^1\}$, and

(c) $\lim_{t\rightarrow\infty}\mu
S(t)=\lambda\mu^1+(1-\lambda)\mu^0$ if and only if
\begin{equation}\label{cond0}
\lim_{t\rightarrow\infty}\sum_y
p_t(x,y){\mu}\{\eta:\eta(y)=1\}=\lambda \text{ for all
}x\in\mathcal{S}.
\end{equation}
\end{thm}

We will say that the transition rates $q(x,y)$ on $\mathbb{Z}^d$
have finite range $N$ if $q(x,y)=0$ when $|x-y|>N$. In order to
show that (\ref{schwartzcond}) is not an unreasonable condition
the following corollary gives circumstances under which
(\ref{schwartzcond}) holds.
\begin{cor}\label{cor3.2}
Let $\mathcal{S}=\mathbb{Z}^d$, $q_e(x,y)=q_e(0,y-x)$, and
$q_v(x,y)=q_v(0,y-x)$. Suppose $X(t)-Y(t)$ is recurrent and
$q_e(x,y)$ has finite range $N$. Then
$\mathcal{I}_e=\{\mu^0,\mu^1\}$ and for $\mu\in\mathcal{P}$,
$\lim_{t\rightarrow\infty}\mu S(t)=\lambda\mu^1+(1-\lambda)\mu^0$
if and only if (\ref{cond0}) holds.
\end{cor}

\begin{thm}\label{thm3.4}
(a) $\mu_\alpha$ exists for all $\alpha\in\mathcal{H}$, and
$\mu_{\alpha_1}=\mu_{\alpha_2}$ if and only if
$\alpha_1R\alpha_2$.

(b) If $g(x,y)<1$ for some $x,y\in \mathcal{S}$ and
\begin{equation}\label{schwartz3} P^E[\int_0^\infty q_e(E_t)\,
dt=\infty]=0 \text{ for some }E\in \mathcal{S}_2
\end{equation}
then $\mathcal{I}_e=\{\mu_\alpha :\alpha\in\mathcal{H}^*_R\}$.

(c) If $q(x,y)=q(y,x)$ for all $x,y\in \mathcal{S}$ and
\begin{equation}\label{schwartz2} P^E[\int_0^\infty q_v(E_t)\,
dt=\infty]=0 \text{ for some }E\in \mathcal{S}_2
\end{equation}
then $\mathcal{I}_e=\{\mu_\alpha :\alpha\in\mathcal{H}_R\}$.
\end{thm}
The condition that $g(x,y)<1$ for some $x,y\in \mathcal{S}$ is not
needed in part (a), but we put it there because if $g\equiv 1$
then we are left with the situation in Theorem \ref{thm3.1}. It
should also be remarked that if $q(x,y)=q(y,x)$ and $g(x,y)<1$ for
some $(x,y)\in \mathcal{S}_2$ then Lemma VIII.1.23 in IPS implies
that (\ref{schwartz3}) and (\ref{schwartz2}) are satisfied.  On
the other hand when $q(x,y)=q(y,x)$, we claim that $g\equiv 1$
implies that $X(t)$ is recurrent so that $\beta(x)+\delta(x)>0$
for some $x$ gives us (\ref{eq0}). To prove the claim use the
Chapman-Kolmogorov equation to get
\begin{eqnarray*}
p_{2t}(x,x)&=&\sum_yp_t(x,y)p_t(y,x)\\
&=&\sum_y[p_t(x,y)]^2=P^{(x,x)}[X_1(t)=X_2(t)].
\end{eqnarray*}
So if $X(t)$ is transient then $g(x,y)<1$ for some $x,y\in
\mathcal{S}$ since
\begin{equation*}
\int_0^\infty P^{(x,x)}[X_1(t)=X_2(t)]\, dt<\infty
\end{equation*}
(This argument will be made more explicit by Lemma
\ref{lemma3.6}).

\begin{thm}\label{thm3.6}
Suppose $\mu\in\mathcal{P}$ and that
$E^{(x,y)}g(X(t),Y(t))\rightarrow 0$ for some $x,y\in
\mathcal{S}$. If
\begin{equation}\label{cond1}
\lim_{t\rightarrow\infty}\sum_y
p_t(x,y){\mu}\{\eta:\eta(y)=1\}=\alpha(x) \text{ and}
\end{equation}
\begin{equation}\label{cond2}
\lim_{t\rightarrow\infty}\sum_{u,v}p_t(x,u)p_t(x,v){\mu}\{\eta:\eta(u)=\eta(v)=1\}=\alpha^2(x)
\text{ for all }x\in\mathcal{S}
\end{equation}
then $\lim_{t\rightarrow\infty} \mu S(t)=\mu_\alpha$. A necessary
and sufficient condition for $\lim_{t\rightarrow\infty} \mu
S(t)=\mu_\alpha$ is that
\begin{equation}\label{cond4}
\lim_{s\rightarrow\infty}\limsup_{t\rightarrow\infty} \int_X
\{\sum_x p_s(w,x)P^x[\Lambda]\sum_y
p_t(x,y)[\eta(y)-\alpha(y)]\}^2\, d\mu(\eta)=0
\end{equation}
\end{thm}
We should mention two instances for which $g(x,y)<1$ for some
$x,y\in\mathcal{S}$ implies
\begin{equation}\label{cond6}
E^{(x,y)}g(X(t),Y(t))\rightarrow 0 \text{ for some }x,y\in
\mathcal{S}.
\end{equation}
Firstly, if $q(x,y)$ is symmetric then as stated in the comments
following Theorem \ref{thm3.4}, Lemma VIII.1.18 in IPS gives
(\ref{cond6}).  Secondly, if the only bounded harmonic functions
are constants then Corollary II.7.3 in IPS together with
Proposition 5.19 in Kemeny, Snell, and Knapp(1976) give
(\ref{cond6}). We also note here that condition (\ref{cond4}) is
equivalent to (\ref{cond1}) and (\ref{cond2}) when
$P^x[\Lambda]=1$ for all $x\in\mathcal{S}$.

\begin{cor}\label{cor1.8}
If $g(x,y)<1$ for some $x,y\in \mathcal{S}$ and
$\mathcal{H}=\{\alpha:\alpha\in[0,1]\}$, then
$\mathcal{I}_e=\{\mu_\alpha :\alpha\in\mathcal{H}\}$.
\end{cor}

The proofs of the above theorems appear in Section \ref{thms}. The
above theorems give partial results concerning the invariant
measures and their respective domains of attraction for certain
NVE processes. Clearly there are NVE processes which are not
covered by these theorems. Examples of these situations include
the process on $\mathbb{Z}^2$ where $q_e(x,y)$ is translation
invariant, $\beta(x)=\delta(x)\equiv 0$, $q_v(x,y)=0$ outside of a
finite set, and $q_v(x,y)$ is not symmetric. A more interesting
example is provided in V.1.6 of IPS; in fact using Liggett's
example we can create similar examples to show that there exist
NVE processes which do not satisfy (\ref{schwartz3}) yet have
$g(x,y)<1$ for some $x,y\in\mathcal{S}$. Section \ref{conjecture}
discusses how one might go about proving a general result that
would include the exceptions we have just mentioned.

We now turn to a discussion of a slightly more general process. In
particular, modify the NVE process by allowing for exclusion rates
where $q_e(x,y)\neq q_e(y,x)$. Call such a process a generalized
NVE process.  It should be noted that not requiring the symmetry
of $q_e(x,y)$ really does change the nature of the process.  We
will state two main reasons for this.  Firstly, the properties of
the dual finite particle system that allow us to prove the above
theorems no longer exist. Secondly, the results for the asymmetric
case are completely different; in fact it is known that Theorems
\ref{thm3.4} and \ref{thm3.6} and Corollary \ref{cor1.8} do not
hold in general when $q_e(x,y)$ is not symmetric. We can however
prove certain things about the generalized NVE process in specific
cases using methods other than duality.

In Section \ref{ergodic} we prove an ergodic theorem for the case
where $q_v(x,y)\equiv 0$ using the coupling method.  When
$q_v(x,y)\equiv 0$ we will call the process a noisy exclusion
process. We will also show in this final section that Theorem
\ref{thm3.0} does not hold in general when $q_e(x,y)$ is not
symmetric.

The main result of Section \ref{ergodic} is an extension, in the
case where $\mathcal{S}=\mathbb{Z}^d$ and the transition rates
have finite range, of Schwartz's(1976) ergodic theorem which is
exactly Theorem \ref{thm3.0} when $q_v(x,y)\equiv 0$. Before we
state the theorem we need the following definitions:
\begin{eqnarray*}
&&T_n=\{x\in \mathbb{Z}^d:|x_i|\le n\text{ for all }i\}.\\
&&T_n^N=T_{n+N}\backslash T_n.
\end{eqnarray*}

\begin{thm}\label{thm2.1}
Suppose $\eta_t$ is a noisy exclusion process with transition
rates $q_e(x,y)$ irreducible with respect to $\mathbb{Z}^d$ and
having finite range $N$. Let $\{b_l\}$ be a nonnegative sequence
satisfying $\text{(a)}$ $\sum b_l=\infty$ if $d=1$ and
$\text{(b)}$~ $\lim_{l\rightarrow\infty}{lb_l}=\infty$ if $d\ge
2$.  If $p(l)$ is a nonnegative function on $\mathbb{N}$
satisfying $p(l+1)\ge p(l)+N$ and is bounded by $kl^k$ for some
$k>0$, and if $\beta,\delta$ satisfy $\beta(x)+\delta(x)\ge b_l$
for all $x\in T_{p(l)}^N$ and $\beta(x)=\delta(x)=0$ otherwise,
then $\eta_t$ is ergodic.
\end{thm}
For some simple examples to see the applicability of Theorem
\ref{thm2.1} set $N=1$ and let $p(l)$ be an arithmetic sequence
e.g. $k,2k,3k,\ldots$ Suppose $\beta(x)=\delta(x)=1$ for all
$\|x\|=nk,\, n\in\mathbb{N}$ with  $\|\cdot\|$ being the
$l^\infty$ norm and $\beta(x)=\delta(x)=0$ otherwise. Then the
theorem tells us that the noisy exclusion process is ergodic. Note
that if $k=1$ and $\beta(x)=\delta(x)=1+\delta$ for $\delta>0$
then the $M<\epsilon$ Theorem in Section I.4 of IPS also gives us
ergodicity for doubly stochastic transition kernels. If $k>1$ then
the $M<\epsilon$ Theorem in general gives us no information. Also,
Theorem \ref{thm2.1} allows us to let
$\beta(x)+\delta(x)\rightarrow 0$ whereas the $M<\epsilon$ Theorem
again gives no information in such a circumstance. We should
however mention here that if $q_e(x,y)$ is symmetric and $k=1$, a
version of the $M<\epsilon$ Theorem proven in Ferarri(1990) allows
for $\beta(x)+\delta(x)\rightarrow 0$, but once again, Ferarri's
theorem gives no information in the case where $k>1$.

\section{The dual process: a finite particle system}
In order to prove the theorems we will need many lemmas. The lemma
in this section which concerns the dual process is the most
important and is in fact the reason that we are able to prove
anything about these processes. Its proof follows that of Theorem
VIII.1.1 in IPS.
Before stating and proving Lemma \ref{lem3.3} we will need some
more definitions.

Let $Y$ be the class of all finite subsets of $\mathcal{S}$
excluding the empty set. The semi-dual process $A_t$ is a
continuous time Markov chain on $Y$ such that the particles in
$A_t$ move independently on $\mathcal{S}$ according to the motions
of the independent $X_i(t)$ processes except that transitions to
sites that are already occupied are handled in the following way:
If a particle at $x$ attempts to move to $y$ which is already
occupied then the transition is either suppressed with probability
$q_e(x,y)/[q_e(x,y)+q_v(x,y)]$ or the two particles coalesce and
move together thereafter with probability
$q_v(x,y)/[q_e(x,y)+q_v(x,y)]$.  In particular $|A_t|\le|A_{t+s}|$
for all $s\ge 0$.

Now let $Y^*$ be defined by adding to $Y$ a cemetery state,
$\Delta$, and the empty set, $\emptyset$. We define the process
$A_t^*$ starting in a state $A\in Y$ to move just as $A_t$ does
except that in addition $A_t^*$ goes to $A_t^*\backslash\{x\}$ at
rate $\beta(x)$ if $x\in A_t^*$ and $A_t^*$ goes to the cemetery
state $\Delta$ at rate $\sum_{x\in A_t^*} \delta(x)$. We will call
$A_t^*$ the dual process. Define $D$ to be the event that $A_t^*$
is never in the state $\Delta$.

If $\mu\in\mathcal{P}$ and $A\in Y$, then define
\begin{equation*}
\hat{\mu}(A)=\mu\{\eta:\eta(x)=1 \text{ for all } x\in A\}.
\end{equation*}
Extend this function to $Y^*$ by letting $\hat{\mu}(\Delta)=0$ and
$\hat{\mu}(\emptyset)=1$.

\begin{lem}\label{lem3.3}
Extend the domain of $\eta\in X$ by letting $\eta(\Delta)=0$. If
$A\in Y$ then for all $t\ge 0$
\begin{equation*}
P^\eta[\{\eta_t=1 \text{ on }A\}]=P^{A}[\{\eta=1 \text{ on
}A_t^*\}\cup\{A_t^*=\emptyset\}].
\end{equation*}
\end{lem}

\begin{proof}
Let
\begin{equation*}
u_\eta (t,A)=P^\eta[\{\eta_t=1 \text{ on
}A\}\cup\{A=\emptyset\}]=S(t)H(\cdot,A)(\eta),
\end{equation*}
where for $A\neq \emptyset$
\begin{equation*}
H(\eta,A)=\prod_{x\in A}\eta(x)=\left\{
\begin{array}{ll}
1&  \text{if }\eta(x)=1 \text{ for all }x\in A^*\\
0&  \text{otherwise}
\end{array}\right.
\end{equation*}
and $H(\eta,\emptyset)=1$.

For each $A\in Y$, $H(\cdot,A)\in\mathcal{D}$ so we have
\begin{eqnarray*}
{\Omega}H(\cdot,A)(\eta)&=&\sum_{\eta(x)=1,\eta(y)=0}{q_e(x,y)[H(\eta_{xy},A)-H(\eta,A)]}\\
&+&\sum_{x,y:\eta(x)\neq
\eta(y)}{q_v(x,y)[H(\eta_{x},A)-H(\eta,A)]}\\
&+&\sum_x
[\beta(x)(1-\eta(x))+\delta(x)\eta(x)][H(\eta_x,A)-H(\eta,A)]\\
&=&\frac{1}{2}\sum_{x,y}{q_e(x,y)[H(\eta_{xy},A)-H(\eta,A)]}\\
&+&
\sum_{x\in A,y\in \mathcal{S}}{q_v(x,y)H(\eta,A\backslash\{x\})[1-2\eta(x)]\{\eta(x)[1-\eta(y)]+\eta(y)[1-\eta(x)]\}}\\
&+& \sum_{x\in A}\beta(x)[H(\eta,A\backslash\{x\})-H(\eta,A)]
+\sum_{x\in A}\delta(x)[H(\eta,\Delta)-H(\eta,A)]
\end{eqnarray*}
\begin{eqnarray*}
&=&\frac{1}{2}\sum_{x,y}{q_e(x,y)[H(\eta,A_{xy})-H(\eta,A)]}+
\sum_{x\in A,y\in \mathcal{S}}{q_v(x,y)H(\eta,A\backslash\{x\})[\eta(y)-\eta(x)]}\\
&+& \sum_{x\in A}\beta(x)[H(\eta,A\backslash\{x\})-H(\eta,A)]
+\sum_{x\in A}\delta(x)[H(\eta,\Delta)-H(\eta,A)]\\
&=&\sum_{x\in A, y\notin A}{q_e(x,y)[H(\eta,A_{xy})-H(\eta,A)]}
+\sum_{x\in A, y\in
\mathcal{S}}q_v(x,y)[H(\eta,(A\backslash\{x\})\cup\{y\})-H(\eta,A)]\\
&+& \sum_{x\in A}\beta(x)[H(\eta,A\backslash\{x\})-H(\eta,A)]
+\sum_{x\in A}\delta(x)[H(\eta,\Delta)-H(\eta,A)].
\end{eqnarray*}
Here $A_{xy}$ is obtained from $A$ in the same way that
$\eta_{xy}$ is obtained from $\eta$.  The symmetry of $q_e(x,y)$
is used in second and fourth steps above.

By Theorem I.2.9 in IPS
\begin{eqnarray*}
&&\frac{d}{dt}u_\eta (t,A)=\Omega S(t)H(\cdot,A)(\eta)\\
&=&\sum_{x\in A, y\notin
A}q_e(x,y)[S(t)H(\cdot,A_{xy})(\eta)-S(t)H(\cdot,A)(\eta)]\\
&+&\sum_{x\in A, y\in
\mathcal{S}}q_v(x,y)[S(t)H(\cdot,(A\backslash\{x\})\cup\{y\})(\eta)-S(t)H(\cdot,A)(\eta)]\\
&+& \sum_{x\in
A}\beta(x)[S(t)H(\cdot,A\backslash\{x\})(\eta)-S(t)H(\cdot,A)(\eta)]
+\sum_{x\in A}\delta(x)[S(t)H(\cdot,\Delta)(\eta)-S(t)H(\cdot,A)(\eta)]\\
 &=&\sum_{x\in A,
y\notin A}q_e(x,y)[u_\eta(t,A_{xy})-u_\eta(t,A)]+\sum_{x\in A,
y\in
\mathcal{S}}q_v(x,y)[u_\eta(t,(A\backslash\{x\})\cup\{y\})-u_\eta(t,A)]\\
&+& \sum_{x\in A}\beta(x)[u_\eta(t,A\backslash\{x\})-u_\eta(t,A)]
+\sum_{x\in A}\delta(x)[u_\eta(t,\Delta)-u_\eta(t,A)].
\end{eqnarray*}
For each $A\in Y$, the unique solution to this system of
differential equations with initial condition $H(\eta,A)$ is
\begin{equation*}
E^AH(\eta,A_t^*)=P^A[\{\eta=1\text{ on
}A_t^*\}\cup\{A_t^*=\emptyset\}]
\end{equation*}
(See Theorem 1.3 of Dynkin(1965)).
\end{proof}

\section{Preliminary lemmas}
The first five lemmas are adaptations of lemmas proved by
Schwartz(1976). We omit the proofs of Lemmas \ref{lemma3.6},
\ref{lemma3.6.5}, and  \ref{lemma3.3.5} since they are the same as
found in Schwartz(1976) except for perhaps a change in notation.

Suppose $\mathcal{E}_t$ is a continuous time nonexplosive jump
process on a countable set $\mathcal{N}$ and let $\mathcal{E}_k$
be the imbedded discrete-time Markov chain. The transition rates
of $\mathcal{E}_t$ are given by $Q_{xy}$. For $\mathcal{L}\subset
\mathcal{N}$ define
\begin{equation*}
Q_\mathcal{L}(x)=\sum_{y\in \mathcal{L},y\neq x} Q_{xy}.
\end{equation*}
\begin{lem}\label{lemma3.6}
Assume there exist constants $0<\alpha_1< \alpha_2<\infty$ such
that for each $x\in \mathcal{N}$, $\alpha_1\le Q_\mathcal{N}(x)\le
\alpha_2$. Then almost surely
\begin{equation*}
\{\omega |\int_0^\infty Q_\mathcal{L}(\mathcal{E}_t)\,
dt=\infty\}=\{\omega | \mathcal{E}_k\in \mathcal{L} \text{
infinitely often}\}\subset\{\omega | \mathcal{E}_t\in \mathcal{L}
\text{ for some t}\}.
\end{equation*}
\end{lem}

\begin{lem}\label{lemma3.6.5}
Assume $0<\sup_x (\beta(x)+\delta(x))<\infty$.  Then (\ref{eq0})
holds if and only if
\begin{equation*}
P^A[A_t^*=\emptyset\text{ or }A_t^*=\Delta\text{ eventually}]=1
\end{equation*}
for all $A\in Y$.
\end{lem}

For the next lemma define the function
\begin{equation*}
h(A)=P^{A}(|A_t|<|A| \text{ for some }t>0)\text{ for }A\in Y
\end{equation*}
which is in some sense a voter model analog of the function
$g(x,y)$.

\begin{lem}\label{lemma3.3.3}
If (\ref{schwartzcond}) holds then $P^{A}(|A_t|=1 \text{
eventually})=1$ for all $A\in Y$.
\end{lem}
\begin{proof}
We first prove the case for which $A_t$ starts in a two particle
state $|A|=2$.

Take $\mathcal{E}_t$ in Lemma \ref{lemma3.6} to be $A_t$, and let
$\mathcal{L}$ be the set of states such that $|A_t|=1$. We then
interpret $Q_\mathcal{L}(A_t)$ as the rate at which $A_t$ jumps to
a one particle state.  If $A_t=\{x\}$ then $Q_\mathcal{L}(A_t)$ is
just $q_x$. Now suppose that $|A_t|=2$ for all $t$. Then $A_t$ is
exactly $E_t$ defined above to be the two particle exclusion
process with respect to $q(x,y)$. Therefore
\begin{equation*}
\int_0^\infty q_v(E_t) \,dt=\int_0^\infty Q_\mathcal{L}(A_t)\,
dt=\infty
\end{equation*}
and by Lemma \ref{lemma3.6}, $|A_t|=1$ eventually, a
contradiction. We have thus proved the case where $|A|=2$.

For the general case suppose $|A|\ge 2$. Couple $B_t$, a semi-dual
process starting from a two particle state $|B|=2$, with $A_t$ so
that $B_t\subset A_t$. In order to do this let $A_t$ and $B_t$
move as usual except when a particle tries to move with rate
$q_e(x,y)$ to an occupied site, instead of the motion being
``excluded", let the two particles switch places. Of course this
is the same motion as before, just a different way of thinking of
it.

Using the coupling we have now that $h(A)=1$ for all $|A|\ge 2$.
Thus with probability one, $|A|$ decreases for all $|A|\ge 2$
which proves the lemma.
\end{proof}

Recall that $D$ is the event where $A_t^*$ is never in the state
$\Delta$.
\begin{lem}\label{lemma3.3.4}
If $\beta(x)\equiv 0$ then
\begin{equation*}
\lim_{t\rightarrow\infty}E^{\{x\}}P^{A_t}[D^c, \Lambda]=0\text{
for all }x\in\mathcal{S}.
\end{equation*}
\end{lem}
\begin{proof}
Let $\mathcal{E}_t=(X(t),\zeta(t))$ be a Markov jump process on
$\mathcal{N}=\mathcal{S}\times\{0,1,2,\ldots\}$ with jump rates
$Q_{(x,n),(y,0)}=q(x,y)$ and $Q_{(x,n),(x,n+1)}=\delta(x)$. Let
$\mathcal{L}=\mathcal{S}\times\{1,2,\ldots\}$ so that
$Q_{\mathcal{L}}((x,n))=\delta(x)$. We then have that
\begin{eqnarray*}
\lim_{t\rightarrow\infty}E^{\{x\}}P^{A_t}[D^c, \Lambda]
&=&\lim_{t\rightarrow\infty}P^{x}[\mathcal{E}_s \text{ jumps to
}\mathcal{L}
\text{ after time }t, \Lambda]\\
&=&P^{x}[\mathcal{E}_k\in\mathcal{L} \text{ infinitely often},
\Lambda].
\end{eqnarray*}
But the right-hand side is equal to $0$ by Lemma \ref{lemma3.6}
completing the proof.
\end{proof}

We will need three definitions in stating the next lemma and in
proving Theorem \ref{thm3.4.1}. Before stating the definitions we
ask the reader to think of $\mu\{\eta:\eta(X(t))=0\}$ as a family
of random variables (indexed by $t$) on the space of paths. We
then have
\begin{equation*}
\mathcal{P}'=\{\mu\in\mathcal{P}:\lim_{t\rightarrow\infty}\mu\{\eta:\eta(X(t))=0\}=1
\text{ almost surely on }\Lambda^c\}.
\end{equation*}
\begin{equation*}
\mathcal{H}'=\{\alpha\in\mathcal{H}:\lim_{t\rightarrow\infty}\alpha(X(t))=0
\text{ almost surely on }\Lambda^c\}.
\end{equation*}
 If $S(t)$ is the
semigroup for an NVE process then let $S'(t)$ be the semigroup for
the same process except that $\beta(x)=\delta(x)\equiv 0$.

For part (b) of the following lemma we couple $A_t$ and $A_t^*$ so
that they move together until the first time that $A_t^*=\Delta$
or $|A_t^*|<|A_t|$.
\begin{lem}\label{lemma3.3.5}
(a) $\mathcal{H}'$ is a set of class representatives for the
equivalence relation $R$ on $\mathcal{H}$.

(b) If we extend the state space of $A_t$ to include $\Delta$ and
$\emptyset$ as absorbing states then
$\lim_{t\rightarrow\infty}P^{{A}_t^*}[{A}_s^*\neq A_s \text{ for
some }s\ge 0]=0$ almost surely.

(c) Suppose that $\beta(x)\equiv 0$. If $\mu\in\mathcal{I}$ or if
$\mu=\lim_{t\rightarrow\infty}\nu S'(t)$ exists for
$\nu\in\mathcal{I}$, then $\mu\in\mathcal{P}'$.
\end{lem}

Define $E_t^n$ to be the finite exclusion process on $n$ particles
starting in the state $A$ where $|A|=n$. To be consistent with our
previous definition of $E_t$ we will leave the superscript off if
$n=2$ so that $E_t=E_t^2$ and $|E|=2$.
\begin{lem}\label{lem3.9}
If (\ref{schwartz2}) holds and $q(x,y)=q(y,x)$  then
\begin{equation*}
P^A[\int_0^\infty \sum_{E\subset E_t^n} q_v(E)\, dt=\infty]=0
\text{ for all }A\in Y.
\end{equation*}
\end{lem}
\begin{proof}
Suppose $A=\{x_1,\ldots,x_n\}$. Let $E^{\{i,j\}}_t$ be the two
particle exclusion process starting from $\{x_i,x_j\}$. We will
show there exists a multiple coupling of the processes $E_t^n$ and
$E^{\{i,j\}}_t$ for $0\le i<j\le n$ such that
\begin{equation}\label{eq3.5}
 \{E_t^n\}\subset\bigcup_{0\le i<j\le n}
\{E^{\{i,j\}}_t\}.
\end{equation}

Let $X_i(t)$ be a process equal in distribution to $X(t)$. The key
to seeing why (\ref{eq3.5}) is true is noticing that there exists
a way to couple $X_i(t)$ and $X_j(t)$ so that whenever one tries
to coalesce with the other, they simply switch places.  This can
be done since $q(x,y)=q(y,x)$. With that said, it is clear that we
can couple the $X_i(t)$'s with $E^n_t$ so that
\begin{equation*}
\{E_t^n\}=\{X_1(t),\ldots,X_n(t)\}.
\end{equation*}
Here the processes $X_i(t)$ start at $x_i$ and are clearly not
independent of each other.

For each $E^{\{i,j\}}_t$ we can label one particle first class and
the other particle second class. We can now think of the evolution
of $E^{\{i,j\}}_t$ in the following way. If a second class
particle tries to go to a site occupied by a first class particle,
it is not allowed to do so.  However, if a first class particle
attempts to move to a site occupied by a second class particle,
the two particles switch places. With this evolution a first class
particle is equal in distribution to $X(t)$.  By choosing the
first class particles to have the paths of the $X_i(t)$ processes
above it is clear that (\ref{eq3.5}) holds.

Suppose now that
\begin{equation*}
P^A[\int_0^\infty \sum_{E\subset E_t^n} q_v(E)\, dt=\infty]>0
\text{ for some }A\in Y.
\end{equation*}
In light of (\ref{eq3.5}), it must be that
\begin{equation*}
P^E[\int_0^\infty q_v(E_t)\, dt=\infty]>0 \text{ for some
}E\in\mathcal{S}_2.
\end{equation*}
By irreducibility
\begin{equation*}
P^E[\int_0^\infty q_v(E_t)\, dt=\infty]>0 \text{ for all
}E\in\mathcal{S}_2.
\end{equation*}
\end{proof}

\begin{lem}\label{lemma3.9}
If $q(x,y)=q(y,x)$ and (\ref{schwartz2}) holds then
$h(E_t^n)\rightarrow 0$ almost surely for all initial states
$A\in\mathcal{S}_n$.
\end{lem}

\begin{proof}
By Lemma \ref{lem3.9},
\begin{equation}\label{eq3.12}
P^A[\int_0^\infty \sum_{E\subset E_t^n} q_v(E)\, dt=\infty]=0
\text{ for all }A\in Y.
\end{equation}

Let $E_k^n$ be the imbedded Markov chain for the process $E_t^n$
starting with initial state $A$. Let $\Omega$ be the path space
for $E_k^n$ and let $\mathcal{M}$ be the probability measure on
$\Omega$ for our process. Choose $\epsilon>0$. If there exists a
set $F\subset\Omega$ such that $\mathcal{M}(F)>0$ and
$h(E_k^n)>\epsilon$ infinitely often on $F$ then it must be that
\begin{equation*}
\sum_{k=0}^\infty \sum_{E\subset E_k^n} q_v(E)=\infty
\end{equation*}
almost surely on $F$ since whenever $\sum_{k=0}^\infty
\sum_{E\subset E_k^n} q_v(E)<\infty$ it must be that
$h(E_k^n)>\epsilon$ finitely many times.

We claim that
\begin{equation}\label{eq3.12.1}
\{\omega|\int_0^\infty \sum_{E\subset E_t^n} q_v(E)\,
dt=\infty\}=\{\omega|\sum_{k=0}^\infty \sum_{E\subset E_k^n}
q_v(E)=\infty\}
\end{equation}
almost surely. To see this define $\tau_k$ to be the $k$th jump
time of $E_t^n$.  Now note that
\begin{equation*}
\int_0^\infty \sum_{E\subset E_t^n} q_v(E)\, dt=\sum_{k=0}^\infty
\sum_{E\subset E_k^n} q_v(E)[\tau_{k+1}-\tau_k].
\end{equation*}
By our assumptions $E[\tau_{k+1}-\tau_{k}]$ and
Var$[\tau_{k+1}-\tau_{k}]$ are bounded above and below uniformly
in $k$.  Since $[\tau_{k+1}-\tau_k|E_1^n, E_2^n, \ldots]$ are
independent, Kolmogorov's Three Series Theorem proves the claim.

Since (\ref{eq3.12.1}) contradicts (\ref{eq3.12}) we have shown
that $h(E_k^n)\rightarrow 0$ almost surely. This however implies
that $h(E_t^n)\rightarrow 0$ almost surely.
\end{proof}

Suppose $V_t$ is the dual process for the voter model with rates
$q(x,y)$ starting from the set $A$. If we couple $A_t$ and $V_t$
so that they move together as much as possible then we can define
the function
\begin{equation*}
f(A)=P^A[A_t\neq V_t\text{ for some }t>0].
\end{equation*}
Again, $f(A)$ plays much the same role as $h(A)$ and $g(x,y)$.
\begin{lem}\label{lemma3.3.6}
If (\ref{schwartz3}) holds then $E^Af(A_t)\rightarrow 0$ for all
$A\in Y$.
\end{lem}

\begin{proof}
We prove first the case where $|A|\le 2$. Let
$\mathcal{E}_t=(A_t,\zeta(t))$ be a Markov jump process on
$\mathcal{N}=(\mathcal{S}_2\cup\mathcal{S})\times\{0,1,2,\ldots\}$
with jump rates (i) $Q_{(A,n),(B,0)}$ equal to the jump rate from
$A$ to $B$ of the semi-dual process
 and (ii) $Q_{(A,n),(A,n+1)}=q_e(A)$ if $|A|=2$. Let
$\mathcal{L}=\mathcal{S}_2\times\{1,2,\ldots\}$ so when $|A|=2$,
$Q_{\mathcal{L}}((A,n))=q_e(A)$ and when $|A|=1$,
$Q_{\mathcal{L}}((A,n))=0$. We then have that
\begin{eqnarray*}
\lim_{t\rightarrow\infty}E^Af(A_t)
&=&\lim_{t\rightarrow\infty}P^A[\mathcal{E}_s \text{ jumps to
}\mathcal{L}
\text{ after time }t]\\
&=&P^{A}[\mathcal{E}_k\in\mathcal{L} \text{ infinitely often}].
\end{eqnarray*}
Since (\ref{schwartz3}) holds, Lemma \ref{lemma3.6} implies that
the right-hand side is $0$.

Now suppose $|A|>2$. Change the coupling of the $X_i(t)$ processes
that we used in Lemma \ref{lem3.9} by letting $X_i(t)$ and
$X_j(t)$ switch places at rate $q_e(X_i(t),X_j(t))$ and coalesce
and move together thereafter at rate $q_v(X_i(t),X_j(t))$. Again,
we are allowed to do this since $q_e(x,y)=q_e(y,x)$. With this new
coupling we can couple the $X_i(t)$'s with $A_t$ so that
\begin{equation*}
\{A_t\}=\{X_1(t),\ldots,X_n(t)\}.
\end{equation*}
As in Lemma \ref{lem3.9}, we use the idea of first class particles
along with the fact that $X_i(t)$ can be coupled with
$E^{\{i,j\}}_t$ so that $\{X_i(t)\}\subset \{E^{\{i,j\}}_t\}$, we
have that the proof for $|A|\le 2$ implies the proof for all $A\in
Y$.
\end{proof}

The next theorem is actually a special case of Theorem
\ref{thm3.4}. We prove this special case right now in order make
the proof of the general case easier to read.
\begin{thm}\label{thm3.4.1}
Suppose $q_e(x,y)\equiv 0$.

(a) $\mu_\alpha$ exists for all $\alpha\in\mathcal{H}$, and
$\mu_{\alpha_1}=\mu_{\alpha_2}$ if and only if
$\alpha_1R\alpha_2$.

(b) $\mathcal{I}_e=\{\mu_\alpha:\alpha\in\mathcal{H}^*_R\}$.
\end{thm}

\begin{proof}
The proof is virtually the same as that of Theorem
1.3 in Schwartz(1976), but it is included here for completeness.  We
will however leave out some repetitive details.

 Let $\mathcal{J}$
represent the set of invariant measures for the case where
$\beta(x)=\delta(x)\equiv 0$, in other words the voter model. In
Chapter V of IPS, it is shown that
$\mathcal{J}_e=\{\mu_\alpha:\alpha\in\mathcal{H}^*\}$. Consider a
certain subset of $\mathcal{J}$, namely
\begin{equation*}
\mathcal{J}'=\{\mu\in\mathcal{J}:
\lim_{t\rightarrow\infty}\mu\{\eta:\eta(X(t))=0\}=1 \text{ almost
surely on }\Lambda^c\}.
\end{equation*}
The main part of the proof is showing that there exists a
bijective affine map between $\mathcal{J}'$ and $\mathcal{I}$. To
avoid confusion, we will put a bar over the extremal invariant
measures of the pure voter model so that we have
$\mathcal{J}_e=\{\bar{\mu}_\alpha:\alpha\in\mathcal{H}^*\}$.

In order to do this we will first consider the case where
$\beta(x)\equiv 0$, but $\delta(x)\ge 0$. We start by coupling the
semi-dual process $A_t$ with $n$ independent processes
$X_1(t),\ldots, X_n(t)$ which start from $A=\{x_1,\ldots,x_n\}$
and are equal in distribution to $X(t)$. In particular, couple the
processes
 so that $A_t\subset\{X_1(t),\ldots,
X_n(t)\}$.  Let $X_i^*(t)$ be the dual process starting from
$\{x_i\}$ and henceforth define $T(t)$ to be the semi-group for
the voter model.

By coupling the processes $A_t^*$ and $A_t$ so that they move
together as much as possible, it is clear that for any measure
$\mu\in\mathcal{P}$ and any $A\in Y$, $S(t)\hat{\mu}(A)\le
T(t)\hat{\mu}(A)$. Thus if $\mu\in\mathcal{I}$ and
$\nu\in\mathcal{J}'$ then $\hat{\mu}(A)\le T(t)\hat{\mu}(A)$ and $
S(t)\hat{\nu}(A)\le \hat{\nu}(A)$. Applying the respective
semigroups once more to both these inequalities gives
$T(s)\hat{\mu}(A)\le T(t+s)\hat{\mu}(A)$ and
$S(t+s)\hat{\nu}(A)\le S(s)\hat{\nu}(A)$ so that
$\lim_{t\rightarrow\infty} \mu T(t)$ and
$\lim_{t\rightarrow\infty} \nu S(t)$ exist by monotonicity and
duality.

Now take  $\mu_1\in\mathcal{J}'$. Let
$\lim_{t\rightarrow\infty}\mu_1S(t)=\mu_2$ and define the map
$\sigma(\mu_1)=\mu_2$.  We will show that $\sigma$ is an affine
bijection from $\mathcal{J}'$ to $\mathcal{I}$.

Since $\mu_1\in\mathcal{P}'$, it follows that
\begin{eqnarray*}
\lim_{t\rightarrow\infty}|T(t)\hat{\mu}_1(A)-S(t)\hat{\mu}_1(A)|&\le&P^A[\bigcup_{1\le
i\le
n}\{X^*_i(t)=\Delta\text{ eventually}, \Lambda\}]\\
&\le& \sum_{i=1}^nP^{\{x_i\}}[D^c, \Lambda].\nonumber
\end{eqnarray*}
By the definition of $\mu_2$ and by the fact that
$\mu_1\in\mathcal{J}$
\begin{equation*}
|\hat{\mu}_1(A)-\hat{\mu}_2(A)| \le \sum_{i=1}^nP^{\{x_i\}}[D^c,
\Lambda].
\end{equation*}
Applying $T(t)$ to both sides of this last inequality and passing
to the limit gives
\begin{equation*}
\lim_{t\rightarrow\infty}|\hat{\mu}_1(A)-T(t)\hat{\mu}_2(A)|
\le\lim_{t\rightarrow\infty}\sum_{i=1}^n E^{\{x_i\}}P^{A_t}[D^c,
\Lambda].
\end{equation*}
Lemma \ref{lemma3.3.4} says that the right-hand side above is
equal to $0$ so that $\lim_{t\rightarrow\infty} \mu_2T(t)=\mu_1$.
This proves that $\sigma$ is injective. If we think of
$X^*(t)=\Delta$ as an absorbing state where $X^*(t)$ continually
jumps to $\Delta$ at exponential rate one then a similar argument
using Lemma \ref{lemma3.3.5} (c) shows $\sigma$ to be surjective.
To see that $\sigma$ is affine note simply that if
$\mu_1,\nu_1\in\mathcal{J}'$ then
\begin{equation*}
\lim_{t\rightarrow\infty}(\lambda\mu_2+(1-\lambda)\nu_2)S(t)=\lambda\mu_1+(1-\lambda)\nu_1.
\end{equation*}

We have thus far shown that there exists an affine bijection
between $\mathcal{J}'$ and $\mathcal{I}$ for the case $\beta\equiv
0$. For the general case we compare the process $\eta_t$ with
birth rates $\beta(x)$ and death rates $\delta(x)$ to a similar
process $\tilde{\eta}_t$ having the same transition rates except
that the death rates are now
$\tilde{\delta}(x)=\beta(x)+\delta(x)$ and the birth rates are
identically $0$. Let the associated dual process, semigroup, and
set of invariant measures for $\tilde{\eta}_t$ be $\tilde{A}_t^*$,
$\tilde{S}(t)$, and $\tilde{\mathcal{I}}$.

Couple the two dual processes so that they make the same
transitions except when a particle in $A_t^*$ dies off due to a
$\beta(x)$ jump, then $\tilde{A}^*_t$ goes to the state $\Delta$.
Since $\tilde{S}(t)\hat{\mu}(A)\le S(t)\hat{\mu}(A)$, we can
repeat the monotonicity arguments used above to show that for
$\nu_1\in\mathcal{I}$ and $\nu_2\in\tilde{\mathcal{I}}$, the
limits $\lim_{t\rightarrow\infty}\nu_1 \tilde{S}(t)=\nu_2$ and
$\lim_{t\rightarrow\infty}\nu_2 {S}(t)=\nu_1$ exist.  If we can
show that
\begin{equation}\label{eq3.1.5.2}
\lim_{t\rightarrow\infty}E^AP^{A_t^*}[A_s^*\neq \tilde{A}_s^*
\text{ for some }s\ge 0]=0\text{ for all }A\in Y
\end{equation}
and similarly that
\begin{equation}\label{eq3.1.5.3}
\lim_{t\rightarrow\infty}E^AP^{\tilde{A}_t^*}[A_s^*\neq
\tilde{A}_s^*\text{ for some }s\ge 0]=0\text{ for all }A\in Y
\end{equation}
then we can also show that the map $\lim_{t\rightarrow\infty}\nu_2
{S}(t)=\nu_1=\tilde{\sigma}(\nu_2)$ is an affine bijection between
$\tilde{\mathcal{I}}$ and $\mathcal{I}$. If we extend the state
space of $A_t$ as in Lemma \ref{lemma3.3.5} (b) then the following
inequalities combined with Lemma \ref{lemma3.3.5} (b) prove
(\ref{eq3.1.5.2}) and (\ref{eq3.1.5.3}):
\begin{equation*}
P^{\tilde{A}_t^*}[A_s^*\neq \tilde{A}_s^* \text{ for some }s\ge
0]\le P^{{A}_t^*}[A_s^*\neq \tilde{A}_s^* \text{ for some }s\ge
0]\le P^{{A}_t^*}[{A}_s^*\ne A_s \text{ for some }s\ge 0].
\end{equation*}

Our desired affine bijection from $\mathcal{J}'$ to $\mathcal{I}$
 is just $\tilde{\sigma}\circ{\sigma}$. We are now ready to prove
 the two parts of the theorem.  We start with part (a).

To prove $\mu_\alpha$ exists we need only show
\begin{equation}\label{eq3.1.5.4}
\lim_{t\rightarrow\infty} \nu_\alpha
S(t)=\lim_{t\rightarrow\infty}\lim_{s\rightarrow\infty}\lim_{r\rightarrow\infty}\nu_\alpha
T(r)\tilde{S}(s) S(t).
\end{equation}
Let $\bar{\mu}_\alpha=\lim_{r\rightarrow\infty} \nu_\alpha T(r)$
and let $\tilde{\bar{\mu}}=\lim_{s\rightarrow\infty} \bar{\mu}
\tilde{S}(s)$.  We have already argued that these limits exist.
Applying $S(t)$ and passing to the limit in the following
inequalities proves (\ref{eq3.1.5.4}).
\begin{eqnarray*}
&&\lim_{t\rightarrow\infty}|S(t)\hat{\nu}_\alpha(A)-S(t)\hat{\tilde{\bar{\mu}}}_\alpha(A)|\\
&\le& \lim_{t\rightarrow\infty}
|S(t)\hat{\nu}_\alpha(A)-\tilde{S}(t)\hat{\bar{{\mu}}}_\alpha(A)|+
\lim_{t\rightarrow\infty}|\hat{\tilde{\bar{\mu}}}_\alpha(A)-S(t)\hat{\tilde{\bar{\mu}}}_\alpha(A)|\\
&\le& \lim_{t\rightarrow\infty}
|S(t)\hat{\nu}_\alpha(A)-T(t)\hat{\nu}_\alpha(A)|+
\lim_{t\rightarrow\infty}|\hat{\bar{\mu}}_\alpha(A)-
\tilde{S}(t)\hat{{\bar{\mu}}}_\alpha(A)| +
\lim_{t\rightarrow\infty}|\hat{\tilde{\bar{\mu}}}_\alpha(A)-
S(t)\hat{\tilde{\bar{\mu}}}_\alpha(A)|\\
&\le& 3P^{{A}}[{A}_s^*\ne A_s \text{ for some }s\ge 0].
\end{eqnarray*}

Suppose now that $\lim_{t\rightarrow\infty}\nu_{\alpha_1}
S(t)=\lim_{t\rightarrow\infty}\nu_{\alpha_2} S(t)$. We have
\begin{eqnarray*}
\hat{\nu}_{\alpha_i}(\{X(s)\})
&=&\lim_{t\rightarrow\infty}E^{X(s)}\hat{\nu}_{\alpha_i}(\{X^*(t)\})\\
&=&P^{X(s)}[X^*(t)=\emptyset\text{ eventually}]+ E^{X(s)}
(\lim_{t\rightarrow\infty}\hat{\nu}_{\alpha_i}(\{X^*(t)\})1_{\{X^*(t)\neq
\emptyset\,\forall\,t, D\}}).\nonumber
\end{eqnarray*}
But since $P^{X(s)}(\{X^*(t)\neq \emptyset\,\forall\,t,
D\})\rightarrow 1$ on $\Lambda_x$ by the arguments given for Lemma
\ref{lemma3.3.4} and since
$E^{X(s)}(\lim_{t\rightarrow\infty}\hat{\nu}_{\alpha_i}(\{X(t)\}))=\alpha_i(X(s))$,
then it follows that $\alpha_1R\alpha_2$.

For the opposite direction if we assume that $\alpha_1R\alpha_2$,
then
\begin{eqnarray*}
\lim_{t\rightarrow\infty}(\nu_{\alpha_1} S(t)-\nu_{\alpha_2}
S(t))(A) &=&
\lim_{t\rightarrow\infty}E^A(\hat{\nu}_{\alpha_1}(A_t^*))
-\lim_{t\rightarrow\infty}E^A(\hat{\nu}_{\alpha_2}(A_t^*))\\
&=& \lim_{t\rightarrow\infty}E^A(\prod_{x\in A_t^*}
{\alpha_1}(x))-\lim_{t\rightarrow\infty}E^A(\prod_{x\in
A_t^*}{\alpha_2}(x))=0.
\end{eqnarray*}

For part (b) it is enough to show that the extreme points of
$\mathcal{J}'$ are
$\{\bar{\mu}_\alpha\in\mathcal{J}:\alpha\in\mathcal{H}^*\cap\mathcal{H}'\}$.
Then applying (\ref{eq3.1.5.4}) along with Lemma \ref{lemma3.3.5}
(a) completes the proof. To prove
$\mathcal{J}'_e=\{\bar{\mu}_\alpha\in\mathcal{J}:\alpha\in\mathcal{H}^*\cap\mathcal{H}'\}$
note that if $\lambda\pi_1+(1-\lambda)\pi_2=\mu\in\mathcal{J}'_e$
for $\pi_1,\pi_2\in\mathcal{J}$ then $\pi_1,\pi_2\in\mathcal{J}'$
and hence $\pi_1=\pi_2=\mu$.  Therefore
$\mu\in\mathcal{J}_e\cap\mathcal{J}'=\{\bar{\mu}_\alpha:\alpha\in\mathcal{H}^*\cap\mathcal{H}'\}$.
On the other hand if $\alpha\in\mathcal{H}^*$ and
$\bar{\mu}_\alpha\in\mathcal{J}'$, then $\bar{\mu}_\alpha$ is an
extreme point of $\mathcal{J}'$.
\end{proof}

\section{Proofs of the theorems}\label{thms}
\begin{proof}[Proof of Theorem \ref{thm3.0}]
Suppose condition (\ref{eq0}) holds.  By Lemma \ref{lem3.3} we
need only show that for any two measures
$\mu_1,\mu_2\in\mathcal{P}$, the limits
$\lim_{t\rightarrow\infty}S(t)\hat{\mu_i}(A)$ exist and are equal
for all $A\in Y$.  But Lemma \ref{lemma3.6.5} implies that
\begin{equation*}
\lim_{t\rightarrow\infty}S(t)\hat{\mu_i}(A)=P^A[A_t^*=\emptyset\text{
eventually}]
\end{equation*}
which is independent of $\mu_i$ proving one direction of the
theorem.

For the opposite direction suppose that (\ref{eq0}) does not hold.
 Lemma \ref{lemma3.6.5} implies that
$P^A[A_t^*=\emptyset\text{ or }A_t^*=\Delta\text{ eventually}]<1$
for some $A\in Y$. Therefore
\begin{equation*}
\lim_{t\rightarrow\infty}S(t)\hat{\delta_1}(A)=P^A[A_t^*=\emptyset\text{
eventually}]+P^A[A_t^*\neq\emptyset\,\forall\,t, D]
\end{equation*}
is not equal to
\begin{equation*}
\lim_{t\rightarrow\infty}S(t)\hat{\delta_0}(A)=P^A[A_t^*=\emptyset\text{
eventually}]
\end{equation*}
for some $A\in Y$ showing that the process is not ergodic.
\end{proof}

\begin{proof}[Proof of Theorem \ref{thm3.1}]
By Lemma \ref{lem3.3}, $\lim_{t\rightarrow\infty} \delta_1S(t)$
exists since
\begin{equation*}
\lim_{t\rightarrow\infty}S(t)\hat{\delta_1}(A)=1-P^A[D^c].
\end{equation*}
Similarly, $\lim_{t\rightarrow\infty} \delta_0S(t)$ exists since
\begin{equation*}
\lim_{t\rightarrow\infty}S(t)\hat{\delta_0}(A)=P^A[A_t=\emptyset\text{
eventually}]
\end{equation*}
completing the proof of part (a).

Consider now part (b). By Lemma \ref{lemma3.3.3} and a coupling
argument it can be seen that if $\lim_{t\rightarrow\infty}E^x
\hat{\mu}(\{X(t)\})$ exists, it is independent of our choice of
$x$. So now using Lemmas \ref{lem3.3} and \ref{lemma3.3.3}
together with the Strong Markov Property, we have that if the
following limits then
\begin{eqnarray}\label{eq3.7}
\lim_{t\rightarrow\infty}S(t)\hat{\mu}(A)&=&\lim_{t\rightarrow\infty}E^{A}\hat{\mu}(A_t^*)\\
&=&P^A[A_t^*=\emptyset\text{ eventually}]+
P^A[\lim_{t\rightarrow\infty}|A_t^*|=1,
D]\lim_{t\rightarrow\infty}E^x \hat{\mu}(\{X(t)\}).\nonumber
\end{eqnarray}
If $\mu\in\mathcal{I}$ then the limits on the left exist. Since we
have assumed that (\ref{eq0}) does not hold then
$P^A[\lim_{t\rightarrow\infty}|A_t^*|=1, D]>0$ so that the last
limit on the right-hand side exists.

Let $\lambda=\lim_{t\rightarrow\infty}E^x \hat{\mu}(\{X(t)\})$ and
consider the invariant measure
$\mu^\lambda=\lambda\mu^1+(1-\lambda)\mu^0$.  We have that for all
$A\in Y$,
\begin{equation*}
\hat{\mu}^\lambda(A)=E^{A}\hat{\mu}^\lambda(A_t^*)=P^A[A_t^*=\emptyset\text{
eventually}]+\lambda P^A[\lim_{t\rightarrow\infty}|A_t^*|=1, D].
\end{equation*}
Since $P^A[A_t^*=\emptyset\text{ eventually}]$ and $
P^A[\lim_{t\rightarrow\infty}|A_t^*|=1, D]$ do not depend on $\mu$
or $\mu^\lambda$, comparing the above equation with (\ref{eq3.7})
gives us $\mu=\mu^\lambda$ showing that every invariant measure is
a mixture of $\mu^1$ and $\mu^0$.  This proves part (b).

Part (c) follows from the above arguments together with the fact
that $\lim_{t\rightarrow\infty} E^x\hat{\mu}(\{X(t)\})$ is
independent of our choice of $x$.
\end{proof}

\begin{proof}[Proof of Corollary \ref{cor3.2}]
We need only show that the recurrence of $Z(t)=X(t)-Y(t)$ implies
(\ref{schwartzcond}).

Let $R$ be the set of all $y\in\mathbb{Z}^d$ such that $|y|\le N$.
By our assumptions we can choose $z\in \mathcal{S}$ so that
$q_v(0,z)>0$. If $E_t=\{x_t,y_t\}$ is the two particle exclusion
process then we will say that $E_t=z$ if $x_t-y_t=z$ and $E_t\in
R$ if $|x_t-y_t|\le N$.

Since $Z(t)$ is recurrent, $Z(t)$ jumps to $0$ infinitely often
and therefore $X(t)$ and $Y(t)$ meet infinitely often.  If there
are infinitely many jumps of $Z(t)$ to $0$ caused by the
$q_v(x,y)$ rates then (\ref{schwartzcond}) automatically holds by
arguments similar to those of Lemma \ref{lemma3.6}. Thus we will
henceforth assume that there are infinitely many jumps of $Z(t)$
to $0$ caused by the $q_e(x,y)$ rates giving us
\begin{equation*}
P^{\{x,y\}}(Z(t)\in R \text{ for some }t>0)=1\text{ for all }
x,y\in\mathcal{S}.
\end{equation*}
By coupling $Z(t)$ and $E_t$ together until the first time that
$X(t)$ and $Y(t)$ meet, we have in fact that
\begin{equation*}
P^{\{x,y\}}(E_t\in R \text{ for some }t>0)=1\text{ for all
}\{x,y\}\in\mathcal{S}_2.
\end{equation*}
If $E_k$ is the embedded discrete-time Markov Chain for $E_t$ then
the above equation implies that
\begin{equation*}
P^{\{x,y\}}(E_k\in R \text{ infinitely often})=1\text{ for all
}\{x,y\}\in\mathcal{S}_2.
\end{equation*}

For a fixed $\bar{t}>0$ let $m=\min_{x\in R
}\{P^x(E_{\bar{t}}=z)\}$. Since our process is irreducible $m>0$
therefore
\begin{equation}\label{eq3.8.5}
P^{\{x,y\}}(E_k=z \text{ infinitely often})=1\text{ for all
}\{x,y\}\in\mathcal{S}_2.
\end{equation}
Now by the same argument given to show (\ref{eq3.12.1}) in the
proof of Lemma \ref{lemma3.9},
\begin{equation*}
\{\omega|\int_0^\infty q_v(E_t)\,
dt=\infty\}=\{\omega|\sum_{k=0}^\infty q_v(E_k)=\infty\}
\end{equation*}
almost surely. By (\ref{eq3.8.5}) we get that (\ref{schwartzcond})
holds as desired.
\end{proof}

\begin{proof}[Proof of Theorem \ref{thm3.4}] For part (a) we prove
only the case where $\beta\equiv 0$; the general case follows the
proof of Theorem \ref{thm3.4.1}.
 Again let $T(t)$ be the semigroup for the voter
model with rates $q(x,y)$. Chapter V in IPS tells us
$\lim_{t\rightarrow\infty} \nu_\alpha T(t)=\bar{\mu}_\alpha$
exists for all $\alpha\in\mathcal{H}$. By coupling the dual of our
process together with the dual of the voter model so that they
move together as much as possible, it is clear that
$S(t)\hat{\bar{\mu}}_\alpha(A)\le
T(t)\hat{\bar{\mu}}_\alpha(A)=\hat{\bar{\mu}}_\alpha(A)$. Applying
$S(s)$ to both sides gives $S(t+s)\hat{\bar{\mu}}_\alpha(A)\le
S(s)\hat{\bar{\mu}}_\alpha(A)$. This monotonicity shows that
$\lim_{t\rightarrow\infty} \bar{\mu}_\alpha S(t)$ exists.  By
previous arguments we also have that $$\lim_{t\rightarrow\infty}
\nu_\alpha
S(t)=\lim_{t\rightarrow\infty}\lim_{s\rightarrow\infty}\nu_\alpha
T(s)S(t).$$

Concerning the rest of the proof we will only prove part (b) since
the proof of part (c) is basically the same as that of (b) except
for replacing the use of Lemma \ref{lemma3.3.6} with Lemma
\ref{lemma3.9}. Just as in the proof of Theorem \ref{thm3.4.1} the
general idea is to show that there exists a bijective, affine map
$\sigma$ between $\mathcal{I}$ and $\mathcal{J}$ where
$\mathcal{J}_e=\{\lim_{t\rightarrow\infty}\nu_\alpha
T(t):\alpha\in\mathcal{H}^*\}$.

For part (b) we will prove only the case where
$\beta(x)=\delta(x)\equiv 0$ so that $A_t^*=A_t$. The general
result follows from the arguments laid out in the proof of Theorem
\ref{thm3.4.1} except for a slight change in the independent
processes $X_1(t),\ldots, X_n(t)$ starting from
$A=\{x_1,\ldots,x_n\}$. For the proof here we must use the
coupling of the $X_i(t)$ processes that we used in the proof of
Lemma \ref{lemma3.3.6} instead of letting them be independent. We
now prove the case $\beta(x)=\delta(x)\equiv 0$.

Take $\mu\in\mathcal{I}$ and  suppose that both $A_t$ and $V_t$
start with initial set $A$. By coupling the two processes so that
$A_t$ contains $V_t$, we see that
\begin{equation*}
|S(t)\hat{\mu}(A)-T(t)\hat{\mu}(A)|\le f(A)=P^A[A_t\neq V_t\text{
for some }t>0].
\end{equation*}
By the invariance of $\mu$
\begin{equation}\label{eq3.13.1}
|\hat{\mu}(A)-T(t)\hat{\mu}(A)|\le f(A)
\end{equation}
so that
\begin{equation*}
|T(s)\hat{\mu}(A)-T({t+s})\hat{\mu}(A)|\le T(s)f(A).
\end{equation*}

By Lemma \ref{lemma3.3.6} and the fact that $S(s)f(A)\rightarrow
0$ implies that $T(s)f(A)\rightarrow 0$, the right-hand side goes
to $0$.  This in turn shows that $\lim_{t\rightarrow\infty}
T(t)\hat{\mu}(A)$ exists. The duality of the voter model which is
a special case of Lemma \ref{lem3.3}, implies that
$\lim_{t\rightarrow\infty}\mu T(t)=\nu$ exists and is invariant
for the voter model with rates $q(x,y)$.

By passing to the limit in (\ref{eq3.13.1})
\begin{equation*}
|\hat{\mu}(A)-\hat{\nu}(A)|\le f(A).
\end{equation*}
Hence Lemma \ref{lemma3.3.6} tells us
$\lim_{t\rightarrow\infty}\nu S(t)=\mu$.

For $\mu\in\mathcal{I}$, if we define
$\sigma(\mu)=\lim_{t\rightarrow\infty}\mu T(t)=\nu$, then the
above arguments have shown that $\sigma$ is injective. A similar
arguments proves $\sigma$ maps onto $\mathcal{J}$. To see that it
is affine note that
\begin{equation*}
\lim_{t\rightarrow\infty}(\lambda\mu_1+(1-\lambda)\mu_2)T(t)=\lambda\nu_1+(1-\lambda)\nu_2.
\end{equation*}
We now conclude the proof of the case $\beta(x)=\delta(x)\equiv 0$
by showing that for
$\bar{\mu}_\alpha=\lim_{t\rightarrow\infty}\nu_\alpha T(t)$,
\begin{equation*}
\lim_{t\rightarrow\infty}\nu_\alpha
S(t)=\lim_{t\rightarrow\infty}\bar{\mu}_\alpha S(t).
\end{equation*}
Applying $S(s)$ to the following inequality and passing to the
limit proves the above equation.
\begin{eqnarray*}
&&\lim_{t\rightarrow\infty}|S(t)\hat{\nu}_\alpha(A)-S(t)\hat{\bar{\mu}}_\alpha(A)|\\
&\le& \lim_{t\rightarrow\infty}
|S(t)\hat{\nu}_\alpha(A)-T(t)\hat{{\nu}}_\alpha(A)|+
\lim_{t\rightarrow\infty}|\hat{\bar{\mu}}_\alpha(A)-S(t)\hat{\bar{\mu}}_\alpha(A)|
\le2f(A).
\end{eqnarray*}
\end{proof}

\begin{proof}[Proof of Theorem \ref{thm3.6}]
Putting $A=\{x_1,\ldots,x_n\}$ let
$W_n(t)\hat{\mu}(A)=E^A\hat{\mu}(\{X_1(t),\ldots,X_n(t)\})$ be the
semigroup for $n$ independent processes. Then the assumptions of
the theorem tell us $W_2(t)g(x,y)\rightarrow 0$ so that

\begin{equation}\label{eq0.5}
P^{\{x,y\}}[X(t)=Y(t)\text{ infinitely often}]=0.
\end{equation}

The proof that (\ref{cond4}) is necessary and sufficient for
$\lim_{t\rightarrow}\mu S(t)=\mu_\alpha$ is proven in Theorem 8.7
in Schwartz(1976). The only thing to note is that the assumption
that $X(t)$ is transient is needed only to show that when
$q(x,y)=q(y,x)$, (\ref{eq0.5}) holds.

The rest of the proof is similar to the proof of Theorems V.1.9 in
IPS. Assume that $\mu$ satisfies (\ref{cond1}) and (\ref{cond2}).
By Lemma \ref{lem3.3} and the definition of $\mu_\alpha$, it
suffices to show that for each $A\in Y$,
\begin{equation}\label{eq3.4.0}
\lim_{t\rightarrow\infty}E^A\hat{\mu}(A_t^*)=\lim_{t\rightarrow\infty}E^A\prod_{x\in
A_t^*} \alpha(x).
\end{equation}
where we make the convention that $\alpha(\Delta)=0$ and
$\prod_{x\in\emptyset}\alpha(x)=1$.

Conditions (\ref{cond1}) and (\ref{cond2}) are equivalent to the
assertion that for each $x\in\mathcal{S}$
\begin{equation*}
\sum_y p_t(x,y)\eta(y)
\end{equation*}
converges in probability to $\alpha(x)$ with respect to $\mu$.
This in turn is equivalent to
\begin{equation}\label{eq3.4.1}
\lim_{t\rightarrow\infty} E^{\{x_1,\ldots,x_n\}}
\hat{\mu}(\{X_1(t),\ldots,X_n(t)\})=\prod_{i=1}^n \alpha(x_i)
\end{equation}
where the $X_i(t)$ are all independent.

Let $\tau_1$ be the first time that either $X_i(t)=X_j(t)$ for
some $1\le i<j\le n$, $A_t^*=\Delta$, or $|A_t^*|$ decreases.
Still putting $A=\{x_1,\ldots,x_n\}$, let $\tau_2$ be the first
time starting from $A_{\tau_1}^*$ that any of the three events
described above occur unless $A_{\tau}^*=\Delta$ in which case we
will let $\tau_2=\infty$. Continuing in this way we can define
$\tau_k$ for all $k\ge 1$.

By (\ref{eq3.4.1}) and the Strong Markov Property, if the limits
below exist then
\begin{eqnarray}\label{eq3.4.2}
\lim_{t\rightarrow\infty} E^A[\hat{\mu}(A_t^*),
\tau_{1}=\infty]&=& \lim_{t\rightarrow\infty}
E^{A}[\hat{\mu}(\{X_1(t),\ldots,X_n(t)\}),\tau_{1}=\infty]\\
&=& \prod_{i=1}^n\alpha(x_i)-\lim_{t\rightarrow\infty}
E^{A}[\hat{\mu}(\{X_1(t),\ldots,X_n(t)\}),\tau_1<\infty]\nonumber\\
&=& \prod_{i=1}^n\alpha(x_i)-E^A\lim_{t\rightarrow\infty}
E^{(X_1(\tau_1),\ldots,X_n(\tau_1))}[\hat{\mu}(\{X_1(t),\ldots,X_n(t)\}),\tau_1<\infty]\nonumber\\
&=& \lim_{t\rightarrow\infty} E^{A}[\prod_{x\in
A_t^*}\alpha(x),\tau_1=\infty]=\lim_{t\rightarrow\infty}
E^{A}[\prod_{x\in A_t}\alpha(x),\tau_1=\infty].\nonumber
\end{eqnarray}
By the convergence theorem for bounded submartingales
$\lim_{t\rightarrow\infty} \prod_{x\in A_t}\alpha(x)$ exists
almost surely so by the Dominated Convergence Theorem the above
limits exist.

Using the Strong Markov Property  once more we can get
\begin{eqnarray}\label{eq3.4.3}
\lim_{t\rightarrow\infty} E^A[\hat{\mu}(A_t^*), \tau_1<\infty] &=&
E^A\lim_{t\rightarrow\infty}
\left(E^{A_{{\tau}_1}^*}[\hat{\mu}(A_t^*),\tau_{1}<\infty,\tau_{2}=\infty]\right)\\
&+& E^A\lim_{t\rightarrow\infty}\left(
E^{A_{\tau_1}^*}[\hat{\mu}(A_t^*),\tau_2<\infty]\right).\nonumber
\end{eqnarray}
But by the argument given for (\ref{eq3.4.2}) the first term on
the right-hand side above equals
\begin{eqnarray*}
&&E^A\lim_{t\rightarrow\infty} E^{A_{\tau_1}^*}[\prod_{x\in
A_t^*}\alpha(x),
\tau_1<\infty,\tau_2=\infty]\\
&=&\lim_{t\rightarrow\infty} E^{A}[\prod_{x\in A_t}\alpha(x),
\tau_1<\infty,\tau_2=\infty]P[A_{\tau_1}^*\neq\emptyset\neq\Delta]+P[A_{\tau_1}^*=\emptyset].
\end{eqnarray*}
  The second term on the right-hand
side of (\ref{eq3.4.3}) equals
\begin{equation*}
E^A\lim_{t\rightarrow\infty}
\left(E^{A_{{\tau}_2}^*}[\hat{\mu}(A_t^*),\tau_{2}<\infty,\tau_{3}=\infty]\right)+
E^A\lim_{t\rightarrow\infty}\left(
E^{A_{\tau_2}^*}[\hat{\mu}(A_t^*),\tau_3<\infty]\right).
\end{equation*}
Since (\ref{eq0.5}) holds we have that $P[\tau_k=\infty\text{ for
some }k]=1$. By repeated use of the arguments above it follows
that (\ref{eq3.4.0}) holds.
\end{proof}

\begin{proof}[Proof of Corollary \ref{cor1.8}]
Again, we prove only the case $\beta(x)=\delta(x)\equiv 0$ so that
$A_t^*=A_t$.  The general case follows from above arguments.

Take $\mu\in\mathcal{I}$ and again let $W_n(t)$ be the semigroup
for $n$ independent random walks
$\vec{X}(t)=(X_1(t),\ldots,X_n(t))$. Couple $A_t$ and $\vec{X}(t)$
so that they move together until the first time that two
coordinates of $\vec{X}(t)$ meet. We then have that
\begin{equation}\label{eq3.19}
|S(t)\hat{\mu}(A)-W_n(t)\hat{\mu}(A)|\le g(A).
\end{equation}
Since $S(t)\hat{\mu}(A)=\hat{\mu}(A)$ then
\begin{equation*}
|W_n(s)\hat{\mu}(A)-W_n({t+s})\hat{\mu}(A)|\le W_n(s)g(A).
\end{equation*}

 Corollary II.7.3 of IPS tells us that $\vec{X}(t)$ has no
nonconstant bounded harmonic functions. By Proposition 5.19 of
Kemeny, Snell, and Knapp(1976) $W_n(t)g(A)\rightarrow 0$ so that
$\lim_{t\rightarrow\infty} W_n(t)\hat{\mu}(A)$ exists and is
harmonic for the random walk $\vec{X}(t)$ on $\mathcal{S}^n$. Such
harmonic functions are constant so we can write
\begin{equation*}
\lim_{t\rightarrow\infty} W_n(t)\hat{\mu}(A)=\alpha_n\text{ for
}|A|=n.
\end{equation*}

The proof of Theorem 2.6 in Liggett(2002) shows that there exists
a random variable $G$ taking values in $[0,1]$ with moment
sequence $\alpha_n$. Since $\alpha_n\le 1$ we can use Carleman's
Condition to show that the random variables $\sum_y
p_t(x,y)\eta(y)$ with respect to the measure $\mu$ converge in
distribution to $G$.

If $\gamma$ is the probability measure on $[0,1]$ for $G$, let
$\mu_\gamma=\int_0^1 \mu_\alpha\, \gamma(d\alpha)$.  Using the
arguments presented in Theorem \ref{thm3.6} we can show that for
each $A\in Y$,
\begin{equation*}
\lim_{t\rightarrow\infty}E^A\hat{\mu}(A_t)=\lim_{t\rightarrow\infty}EG^{|A_t|}=\lim_{t\rightarrow\infty}E^A\hat{\mu}_\gamma(A_t).
\end{equation*}
Thus $\mu=\mu_\gamma$ and is hence a mixture of the measures
$\{\mu_\alpha:\alpha\in[0,1]\}$. By Theorem \ref{thm3.6}, each
measure $\mu_\alpha$ has a different domain of attraction proving
that $\mathcal{I}_e=\{\mu_\alpha:\alpha\in[0,1]\}$.
\end{proof}

\section{Further results}\label{conjecture}
The brief discussion below shows how one might adapt
Schwartz(1976) and Chapter V in IPS in order to obtain a general
result. Let
\begin{equation*}\label{def3.2}
{\hat{\mathcal{E}}}=\{\omega:\int_0^\infty q_v(E_t)\, dt=\infty
\}.
\end{equation*}
In the introduction we argued that
$\lim_{t\rightarrow\infty}\alpha(X(t))$ exists almost surely so we
can define $\hat{\mathcal{H}}$ to be the set of all
$\alpha\in\mathcal{H}$ such that
\begin{equation*}\label{def3.3}
\lim_{t\rightarrow\infty} \alpha(X(t))=0\text{ or }1\text{ a.s. on
}{\hat{\mathcal{E}}}
\end{equation*}
where $X(t)$ starts from $x$ if $E_0=\{x,y\}$.
 For those that are keeping track,
${\hat{\mathcal{E}}}$ and $\hat{\mathcal{H}}$ are analogous to
$\mathcal{E}$ and $\mathcal{H}^*$.

Following Schwartz(1976) and Chapter V in IPS, we conjecture that
$\mathcal{I}_e=\{\mu_\alpha: \alpha\in\hat{\mathcal{H}}_R\}$. In
order to prove this one would have to generalize Theorem
\ref{thm3.6} and show that for $\mu\in\mathcal{I}_e$
\begin{equation}\label{eq3.11}
\hat{\mu}(\{X(t),Y(t)\})\rightarrow
\hat{\mu}(\{X(t)\})\hat{\mu}(\{Y(t)\}).
\end{equation}

 As mentioned in the
introduction, it is the monotonicity of $S(t)\hat{\nu}_\alpha(A)$
that allows us to do this for the pure voter model or the pure
symmetric exclusion process. If one were to prove (\ref{eq3.11})
and Theorem \ref{thm3.6} in general, new techniques would be
needed.

\section{An ergodic theorem for a related process}\label{ergodic}
The proof of Theorem \ref{thm2.1} requires the following lemma:
\begin{lem}\label{lem2.1}
Suppose $\{a_n\}$ is bounded above by $k_1n^{k_2-1}$ for some
$k_1, k_2>0$ and that $a_n>0$ for all $n$.  Then there exists a
sequence $\{w_n\}$ such that
\begin{equation*}
\text{(i) } \underset{n\rightarrow\infty}{\liminf}
\frac{a_n}{w_n}=1 \text{ and (ii)
}\underset{n\rightarrow\infty}{\limsup}
\frac{nw_n}{\sum_{l=0}^{n-1}{w_l}}<\infty.
\end{equation*}
\end{lem}

\begin{proof}
If for some sequence $\{w_n\}$ we have that $w_l/l^{k_2}\ge
w_n/n^{k_2}$ for $l\le n$ then
\begin{equation*}
\sum_{l=0}^{n-1}w_l\ge\sum_{l=0}^{n-1}{\frac{w_n}{n^{k_2}}l^{k_2}}\ge
\frac{w_n}{n^{k_2}}\frac{(n-1)^{k_2+1}}{k_2+1}
\end{equation*}
so that condition (ii) holds.  So it remains to find a sequence
$\{w_n\}$ satisfying condition (i) and the inequality
$w_l/l^{k_2}\ge w_n/n^{k_2}$ for $l\le n$.  Let $w_0=a_0$ and let
$w_n=w_{n-1}$ unless $a_n/n^{k_2}=\min_{l\le n}{a_l/l^{k_2}}$ in
which case we let $w_n=a_n$. Then $w_l/l^{k_2}\ge w_n/n^{k_2}$ for
$l\le n$. Now since $\{a_n\}$ is bounded above by $k_1n^{k_2-1}$
it follows that $a_n/n^{k_2}\rightarrow 0$ and hence
$a_n/n^{k_2}=\min_{l\le n}{a_l/l^{k_2}}$ infinitely often so that
$w_n=a_n$ infinitely often.  Therefore (i) is also satisfied by
this choice of $\{w_n\}$.
\end{proof}

Using the basic coupling for the exclusion process combined with
the basic coupling for spin systems, we have that the basic
coupling for a noisy exclusion process has generator
\begin{eqnarray*}
\bar{\Omega}f(\eta,\xi)&=&\sum_{\scriptsize{
\begin{array}{c}
\eta(x)=\xi(x)=1\\
\eta(y)=\xi(y)=0
\end{array}}
}{q_e(x,y)[f(\eta_{xy},\xi_{xy})-f(\eta,\xi)]}\\
&+&\sum_{ \scriptsize{
\begin{array}{c}
\eta(x)=1,\eta(y)=0 \text{ and }\\
\xi(y)=1 \text{ or } \xi(x)=0
\end{array}}
}{q_e(x,y)[f(\eta_{xy},\xi)-f(\eta,\xi)]}\\
&+&\sum_{ \scriptsize{
\begin{array}{c}
\xi(x)=1,\xi(y)=0 \text{ and }\\
\eta(y)=1 \text{ or } \eta(x)=0
\end{array}}
}{q_e(x,y)[f(\eta,\xi_{xy})-f(\eta,\xi)]}\\
&+&\sum_{x:\eta(x)\neq\xi(x)}{c_1(x, \eta)[f(\eta_x,\xi)-f(\eta,\xi)]}
+\sum_{x:\eta(x)\neq\xi(x)}{c_2(x, \xi)[f(\eta,\xi_x)-f(\eta,\xi)]}\\
&+&\sum_{x:\eta(x)=\xi(x)}{c(x, \eta,
\xi)[f(\eta_x,\xi_x)-f(\eta,\xi)]}
\end{eqnarray*}
where
\begin{equation*}
\begin{array}{lr}
c_1(x, \eta)=\left\{
\begin{array}{ll}
   \beta(x) &\text{when }\eta(x)=0\\
   \delta(x) &\text{when } \eta(x)=1
\end{array}\right.&
c_2(x, \xi)=\left\{
\begin{array}{ll}
   \beta(x) &\text{when }\xi(x)=0\\
   \delta(x) &\text{when } \xi(x)=1
\end{array}\right.
\end{array}
\end{equation*}
\begin{equation*}
\text {and }c(x, \eta, \xi)=\left\{
\begin{array}{ll}
   \beta(x)&\text{when }\eta(x)=\xi(x)=0\\
   \delta(x)&\text{when } \eta(x)=\xi(x)=1.
\end{array}\right.
\end{equation*}
Let $\mathcal{\bar{I}}$ be the set of invariant measures for this
coupling.

In order to simplify the notation we define the functions
\begin{equation*}
\begin{array}{lll}
 f_x(\eta,\xi)=[1-\eta(x)]\xi(x),&&h_{yx}(\eta,\xi)=[1-\eta(y)][1-\xi(y)]f_x(\eta,\xi),\\
g_{yx}(\eta,\xi)=\eta(y)\xi(y)
f_x(\eta,\xi),&\text{and}&f_{yx}(\eta,\xi)=\eta(y)[1-\xi(y)]f_x(\eta,\xi).
\end{array}
\end{equation*}
In particular, for $T$ a finite subset of $\mathcal{S}$ we have
\begin{eqnarray} \label{eq1}
\tilde{\Omega}\left(\sum_{x\in T} {f_x(\eta,\xi)}\right)
&=&-\sum_{x\in T,y\in\mathcal{S}}
{\left(q_e(x,y)+q_e(y,x)\right)f_{yx}(\eta,\xi)}-\sum_{x\in T}
{\left(\beta(x)+\delta(x)\right)f_x(\eta,\xi)}\\
&+&\sum_{x \in T, y \notin T} {[q_e(x,y)g_{xy}-q_e(y,x)g_{yx}]}
+\sum_{x\in T, y\notin T} {[q_e(y,x)h_{xy}
-q_e(x,y)h_{yx}]}\nonumber.
\end{eqnarray}

\begin{proof}{Proof of Theorem \ref{thm2.1}}
Recall that $T_n=\{x\in \mathbb{Z}^d:|x_i|\le n \}$. Couple two
noisy exclusion processes, $\eta_t$ and $\xi_t$, with
$\nu\in\mathcal{\bar{I}}$ so that
\begin{equation*}
\int{\bar{\Omega}\left(\sum_{x\in
T_n}{f_x(\eta_t,\xi_t)}\right)d\nu}=0.
\end{equation*}
If we let $\int f_x(\eta_t,\xi_t) d\nu=a(x)$ then since
$f_{yx}(\eta_t,\xi_t)\ge 0$, equation (\ref{eq1}) gives us
\begin{eqnarray} \label{eq2}
&&\sum_{x\in T_{n}}{\left(\beta(x)+\delta(x)\right)a(x)}\\
&\le&\sum_{x \in T_{n}, y \notin T_n}{q_e(x,y)\int (g_{xy} -
h_{yx})d\nu}
+\sum_{x\in T_{n}, y\notin{T}_n}{q_e(y,x)\int (h_{xy} - g_{yx}) d\nu}\nonumber\\
&\le&\sum_{x \in T_{n}, y \notin T_n}{q_e(x,y)a(y)} + \sum_{x \in
T_{n}, y \notin T_n}{q_e(y,x)a(y)}\nonumber\\
&\le& \frac{(2N+1)^d}{2}\sum_{y\in T_{n}^N}{a(y)}+\sum_{y\in
T_{n}^N}{a(y)}\le C_1\sum_{y\in T_{n}^N}{a(y)}\nonumber
\end{eqnarray}
for some constant $C_1$.  If we define
\begin{equation*}
{a_l}=\sum_{y\in T_{p(l)}^N}{a(y)}
\end{equation*}
then by the inequality $\beta(x)+\delta(x)\ge b_l$ for $x\in
T_{p(l)}^N$ we can rewrite (\ref{eq2}) as
\begin{equation} \label{eq3}
\sum_{l=0}^{n-1}{b_{l}a_{l}}\le C_1a_{n}.
\end{equation}

Now suppose $d=1$ and condition (a) in the theorem holds.  Then
since $a(x)\le 1$, we have $a_n\le 2N$ for all $n$.  In light of
equation (\ref{eq3}) we then have that $\sum_{l\ge
0}{b_{l}a_{l}}<\infty$. On the other hand, if we multiply both
sides of (\ref{eq3}) by $b_{n}$ and then sum over $n$  we get
\begin{equation*}
\sum_{n\ge 0}{b_{n}}\sum_{l=0}^n{b_{l}a_{l}}\le C_1\sum_{n\ge
0}{b_{n}a_{n}}<\infty.
\end{equation*}
Rewriting the left hand side we get
\begin{equation*}
\sum_{n\ge 0}{b_{n}}\sum_{l=0}^n{b_{l}a_{l}}=\sum_{l\ge
0}{b_{l}a_{l}}\sum_{n\ge l}{b_{n}}< \infty.
\end{equation*}
This implies that $b_{l}a_{l}=0$ for all $l$ since condition (a)
gives us $\sum b_l=\infty$.  So we have $a(x)=\int f_x d\nu=0$ for
all $x$ so that the marginals of $\nu$ are exactly the same.

Suppose now that $d\ge 2$ and that condition (b) of the theorem
holds. Since $p(l)\le kl^k$ we have that $a_{n}$ is bounded above
by $k_1n^{k_2-1}$ for some $k_1, k_2$.  If we assume that for all
$n$, $a_{n}>0$ then by Lemma \ref{lem2.1}, there exists a sequence
$w_n$ such that $\liminf a_{n}/w_n=1$ and $\limsup
nw_n/\sum_{l=0}^{n-1}{w_l}<~\infty$.  By condition (b), we have
then that
\begin{equation}\label{eq2.4}
\liminf nb_na_{n}/w_n=\infty
\end{equation}
However, we also have that there exists a subsequence $\{n_j\}$
for which
\begin{equation}\label{eq2.5}
\sum_{l=0}^{n_j-1}{b_{l}a_{l}}\le C_1a_{n_j}\le C_2w_{n_j}\le
C_3\frac{\sum_{l=1}^{n_j-1}{w_l}}{n_j}\le
C_3\sum_{l=1}^{n_j-1}{\frac{w_l}{l}}.
\end{equation}
Notice now that if the limit of the right hand side is infinite,
(\ref{eq2.4}) and (\ref{eq2.5}) contradict each other so that we
must have $a_{n}=0$ for some $n$ and consequently $a(x)=\int f_x
d\nu=0$ for all $x$ by irreducibility.  If the right hand side is
bounded then we can use the argument given above for the case
$d=1$ to show that $a(x)=\int f_x d\nu=0$ for all $x$.  In either
case we have that the marginals of $\nu$ are the same, and we thus
have ergodicity of the process.
\end{proof}

We now restrict ourselves to the case where $d=1$ and the
transition rates are $q_e(x,x+1)=p>1/2$ and $q_e(x,x-1)=1-p=q<1/2$
for all $x$. In order to show the importance of the condition that
there exist a sequence $b_l$ satisfying $b_l\le
\beta(x)+\delta(x)$ for all $x\in T_{p(l)}^N$, we will find
examples of processes on $\mathbb{Z}$ that are not ergodic but
satisfy $\sum_x \beta(x)=\infty$.

To start off, consider the case where we have $\beta>0$ and
$\delta>0$ for a single fixed $z$ and no births and deaths at any
other site.  Choose $c$ so that
$c\pi(z)/(1+c\pi(z))=\beta/(\beta+\delta)$ for a reversible
measure $\pi(x)$ on $\mathbb{Z}$. The product measure $\nu^c$ with
marginals $\nu^c\{\eta:\eta(x)=1\}=c\pi(z)/(1+c\pi(z))$ is
reversible with respect to the exclusion process, and its marginal
measure at the site $z$ is reversible with respect to the birth
and death process so that $\nu^c$ is reversible with respect to
the noisy exclusion process. The product measure $\nu_{\rho}$ with
marginals $\nu_\rho\{\eta:\eta(x)=1\}=\rho$ where
$\rho=\frac{\beta}{\beta+\delta}$, is also invariant with respect
to the exclusion process, and again, its marginal measure at the
site $z$ is reversible with respect to the birth and death
process.  So $\nu_\rho$ is also invariant with respect to the
noisy exclusion process.

We have two more invariant measures by starting the process off
with initial states $\delta_0$ and $\delta_1$. This is because
some subsequence of $\lim_{n\rightarrow\infty}\frac{1}{T_n}
\int_0^{T_n}\delta_1S(t)dt$ for $T_n\rightarrow\infty$ must lie
above both of the invariant measures we have constructed above.
Similarly some subsequence of
$\lim_{n\rightarrow\infty}\frac{1}{T_n}
\int_0^{T_n}\delta_0S(t)dt$ lies below the two invariant measures.

In order to show that the noisy exclusion process with
$\beta(z_i)>0$ if and only if $\delta(z_i)>0$ for a finite number
of sites $\{z_1,...,z_k\}$ is not ergodic (this is a special case
of Proposition \ref{prop2.1} below) we will need the following
coupling for two noisy exclusion processes with the same
transition and death rates, but different birth rates. If
$\beta_1(x)$ for the process $\eta_t$ is greater than $\beta_2(x)$
for the process $\xi_t$ for all $x$ then we can couple the two
processes in such a way that $\eta_t\ge \xi_t$. Formally, we have
the coupling given by
\begin{eqnarray*}
\bar{\Omega}f(\eta,\xi)&=&\sum_{ \scriptsize{
\begin{array}{c}
\eta(x)=\xi(x)=1\\
\eta(y)=\xi(y)=0
\end{array}}
}{q_e(x,y)[f(\eta_{xy},\xi_{xy})-f(\eta,\xi)]}\\
&+&\sum_{ \scriptsize{
\begin{array}{c}
\eta(x)=1,\eta(y)=0 \text{ and }\\
\xi(y)=1 \text{ or } \xi(x)=0
\end{array}}
}{q_e(x,y)[f(\eta_{xy},\xi)-f(\eta,\xi)]}\\
&+&\sum_{ \scriptsize{
\begin{array}{c}
\xi(x)=1,\xi(y)=0 \text{ and }\\
\eta(y)=1 \text{ or } \eta(x)=0
\end{array}}
}{q_e(x,y)[f(\eta,\xi_{xy})-f(\eta,\xi)]}\\
&+&\sum_{x:\eta(x)\neq\xi(x)}{c_1(x, \eta)[f(\eta_x,\xi)-f(\eta,\xi)]}+\sum_{x:\eta(x)\neq\xi(x)}{c_2(x, \xi)[f(\eta,\xi_x)-f(\eta,\xi)]}\\
&+&\sum_{x:\eta(x)=\xi(x)}{c(x, \eta,
\xi)[f(\eta_x,\xi_x)-f(\eta,\xi)]}+\sum_{x:\eta(x)=\xi(x)=0}{(\beta_1(x)-\beta_2(x))[f(\eta_x,\xi)-f(\eta,\xi)]}
\end{eqnarray*}
where
\begin{equation*}
\begin{array}{cc} c_1(x, \eta)=\left\{
\begin{array}{ll}
   \beta_1(x) &\text{when }\eta(x)=0\\
   \delta(x) &\text{when } \eta(x)=1
\end{array}\right.&
c_2(x, \xi)=\left\{
\begin{array}{ll}
   \beta_2(x) &\text{when }\xi(x)=0\\
   \delta(x) &\text{when } \xi(x)=1
\end{array}\right.
\end{array}
\end{equation*}
\begin{equation*}
\text{and }c(x, \eta, \xi)=\left\{
\begin{array}{ll}
   \beta_2(x) &\text{when }\eta(x)=\xi(x)=0\\
   \delta(x) &\text{when } \eta(x)=\xi(x)=1.
\end{array}\right.
\end{equation*}
Similarly, we can couple two processes together so that
$\eta_t\le\xi_t$ when $\eta_t$ and $\xi_t$ have the same
transition and birth rates, but death rates such that
$\delta_1(x)\ge\delta_2(x)$ for all $x$.

\begin{prop}\label{prop2.1}
Suppose that $q_e(x,x+1)=p>\frac{1}{2}$ and $q_e(x,x-1)=1-p=q$ for
all $x$ and that $\beta(x)>0$ if and only if $\delta(x)>0$. If
there exists a $z$ such that $\beta(x)=0$ for either all $x\le z$
or for all $x\ge z$ and if there exist $a_1$ and $a_2$ such that
$\frac{a_1\pi(x)}{1+a_1\pi(x)}\le
\frac{\beta(x)}{\beta(x)+\delta(x)}\le
\frac{a_2\pi(x)}{1+a_2\pi(x)}$ for all $x$ where $\beta(x)>0$,
then the process is not ergodic.
\end{prop}

\begin{proof}
Without loss of generality suppose that $\beta(x)=0$ for all
positive $x$ and let $\{z_i\}$ denote the set of points where
$\beta(x)>0$. If $\eta_t$ is the process described in the
hypothesis of the proposition, let the process $\xi_t$ be the same
as $\eta_t$ except that we change the death rates of $\xi_t$ so
that
$\frac{\beta(z_i)}{\beta(z_i)+\delta(z_i)}=\frac{a_1\pi(z_i)}{1+a_1\pi(z_i)}$
for all $\{z_i\}$. Let the process $\zeta_t$ be the same as
$\eta_t$ except that we change the birth rates of $\zeta_t$ so
that
$\frac{\beta(z_i)}{\beta(z_i)+\delta(z_i)}=\frac{a_2\pi(z_i)}{1+a_2\pi(z_i)}$
for all $\{z_i\}$. We can triple couple $\xi_t, \eta_t,$ and
$\zeta_t$ so that $\xi_t\le \eta_t\le \zeta_t$.  Since the measure
$\nu^{a_1}$ is invariant for $\xi_t$ and $\nu^{a_2}$ is invariant
for $\zeta_t$, then $\eta_t$ has an invariant measure $\mu_1$ with
$\nu^{a_1}\le \mu_1\le \nu^{a_2}$.

Let
$M=\max_i\left(\frac{\beta(z_i)}{\beta(z_i)+\delta(z_i)}\right)$.
Note that this maximum is achieved since we assumed earlier that
$\beta(x)=0$ for all positive $x$ and consequently if there exist
an infinite number of $z_i$'s then
$\lim_{i\rightarrow\infty}{\frac{\beta(z_i)}{\beta(z_i)+\delta(z_i)}}=0$.
Now let the process $\zeta_t$ be the same as $\eta_t$ except that
we change the birth rates of $\zeta_t$ so that
$\frac{\beta(z_i)}{\beta(z_i)+\delta(z_i)}=M$ for all $\{z_i\}$.
Again, we can couple $\eta_t$ and $\zeta_t$ so that $\eta_t\le
\zeta_t$. The measure $\nu_M$ is invariant for $\zeta_t$. So
$\eta_t$ has an invariant measure $\mu_2$ such that $\mu_2\le
\nu_M$. Since $\mu_2$ is different from $\mu_1$, the process is
not ergodic.
\end{proof}

Note that using the above proposition, we can construct examples
of nonergodic processes that satisfy all of the hypotheses for
Schwartz's ergodic theorem except for $q_e(x,y)=q_e(y,x)$.

\textbf{Acknowledgement}. The author thanks his advisor, Thomas M.
Liggett, for suggesting the problem and for his encouragement
during the undertaking.

\end{document}